\documentclass[12pt]{article}
\usepackage[margin=1in]{geometry}

\usepackage{cite}
\usepackage{amsmath,amssymb,amsfonts}
\usepackage{algorithmic}
\usepackage{caption}
\usepackage{graphicx}
\usepackage{textcomp}
\usepackage{xcolor}
\def\BibTeX{{\rm B\kern-.05em{\sc i\kern-.025em b}\kern-.08em
    T\kern-.1667em\lower.7ex\hbox{E}\kern-.125emX}}

\title{\LARGE \bf
Securing Signal-free Intersections against Strategic Jamming Attacks: A Macroscopic Approach\\}

\author{Yumeng Bai, Saurabh Amin, Xudong Wang, and Li Jin
\thanks{This work was in part supported by UM-SJTU Joint Institute, J. Wu \& J. Sun Endowment Fund, NSFC Project 62103260, and US AFOSR FA9550-19-1-0263.}
\thanks{Y. Bai, X. Wang and L. Jin are with the UM Joint Institute, Shanghai Jiao Tong University, China. Y. Bai is also with the School of Aeronautics and Astronautics, Purdue University, USA. S. Amin is with the Department of Civil \& Environmental Engineering and with the Laboratory for Information \& Decision Systems, Massachusetts Institute of Technology, USA (emails: bai145@purdue.edu, amins@mit.edu, wxudong@sjtu.edu.cn, li.jin@sjtu.edu.cn).}%
}

\begin{document}

\maketitle
\thispagestyle{empty}
\pagestyle{empty}

\begin{abstract}
We consider the security-by-design of a signal-free intersection for connected and autonomous vehicles in the face of strategic jamming attacks. We use a fluid model to characterize macroscopic traffic flow through the intersection, where the saturation rate is derived from a vehicle coordination algorithm. We model jamming attacks as sudden increase in communication latency induced on vehicle-to-infrastructure connectivity; such latency triggers the safety mode for vehicle coordination and thus reduces the intersection saturation rate. A strategic attacker selects the attacking rate, while a system operator selects key design parameters, either the saturation rate or the recovery rate. Both players' actions induce technological costs and jointly determine the mean travel delay. By analyzing the equilibrium of the security game, we study the preferable level of investment in the intersection's nominal discharging capability or recovery capability. 
\end{abstract}

\section{Introduction}
Signal-free intersections are an emerging traffic operation enabled by vehicle-to-infrastructure (V2I) connectivity (e.g., dedicated short-range communication (DSRC)~\cite{DSRC},~\cite{V2I}).
The replacement of traditional signalized intersections with signal-free intersections has the potential to greatly improve key transportation performance metrics, including throughput, travel delay, safety, etc~\cite{intersection},~\cite{signal-free}.
However, such benefits heavily depend on the quality of V2I connectivity, which may be vulnerable to malicious attacks~\cite{attacks},~\cite{vulnerable},~\cite{vanet}. To the best of the authors' knowledge, despite the increasing interest in this field, very limited tools and methods have been developed to analyze the security risk behind V2I-based signal-free intersections. 

\begin{figure}[htbp]
    \centerline{\includegraphics[width=0.5\linewidth]{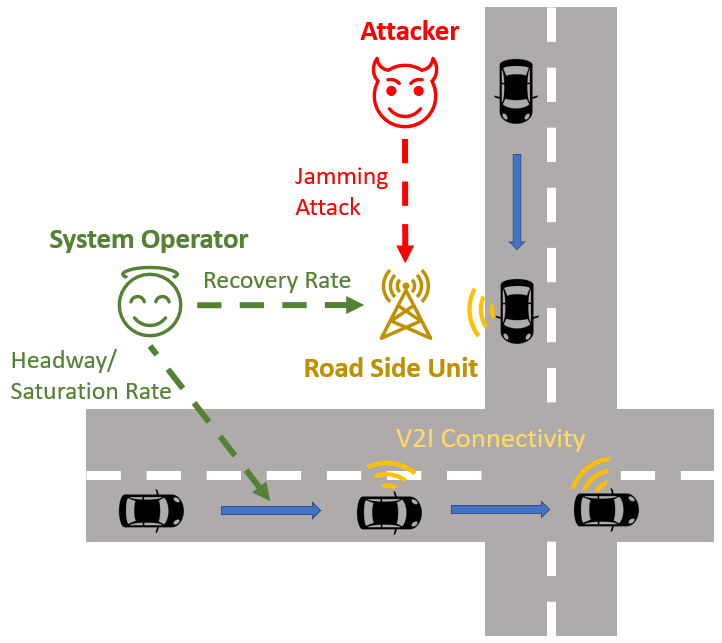}}
    \caption{A typical signal-free intersection.}
    \label{fig:intersec_fig}
\end{figure}

In this paper, we respond to the above challenge by studying the security risk of a representative signal-free intersection configuration based on V2I connectivity, as shown in Fig.~\ref{fig:intersec_fig}. 
In such a system, jamming attacks are one of the major threats that can compromise the flow of data and cause severe communication latency.
The key of our approach is to use a piecewise-deterministic queuing (PDQ) model, which is a macroscopic fluid model, to bridge the gap between lower-level vehicle coordination algorithms and upper-level security risks analysis.
Specifically, the PDQ model captures the jamming attack-induced loss of traffic flow saturation rate and provides analytical utility functions for the security game between a strategic attacker and a strategic system operator (SO).
Analysis of the game's equilibrium provides insights about the security risks in different scenarios and helps determine the optimal level of investment in vehicle coordination hardware/infrastructure and recovery capabilities against jamming attacks.

There has been an increasing interest in security vulnerabilities for connected and autonomous vehicle (CAV) coordination~\cite{cdc18},~\cite{platoongame},~\cite{survey},~\cite{survey2},~\cite{drl},~\cite{drl2}. 
In particular, jamming attack is one of the main forms of Denial-of-Service (DOS) attacks in DSRC. Jamming attacks aim to block legitimate communications and therefore degrade overall system performance, usually in the form of broadcasting radio signal, switching resource blocks, etc~\cite{DOS},~\cite{jamdetect}.
Multiple studies have analyzed different forms of jamming attacks, revealing the severe impact of jamming attacks on communication latency and the difficulty of detection~\cite{DOS},~\cite{jamdetect},~\cite{jamexp},~\cite{jamvul}. However, in the context of signal-free intersections, one outstanding question is: what is the impact of jamming attacks on the intersection performance metrics such as capacity and delay. 
A major challenge for answering the above question is the lack of appropriate model to quantitatively map the attacker's and the SO's behavior to the intersection's performance.

So far, the majority of previous work on signal-free intersections is based on microscopic, trajectory-based models~\cite{latency},~\cite{trajectory},~\cite{trajectory2},~\cite{trajectory3}; very limited work has been done at the macroscopic, intersection level. 
Some researchers have studied potential security risks for CAV against jamming attacks in microscopic settings, such as sensor anomaly detection, access technology, etc~\cite{detect},~\cite{cognitive},~\cite{DRLCAV}. 
However, microscopic models such as trajectory-based ones do not lead to macroscopic performance metrics; an appropriate modeling approach is needed to bridge this gap.
In this regard, macroscopic traffic models have been developed and applied to various CAV-related applications in recent years~\cite{PDQ},~\cite{piecewise},~\cite{piecewise2}.
Nevertheless, to the best of our knowledge, macroscopic security risk analysis for CAV in signal-free intersections facing jamming attacks has not been studied yet.

\begin{figure}[htbp]
    \centerline{\includegraphics[width=0.6\linewidth]{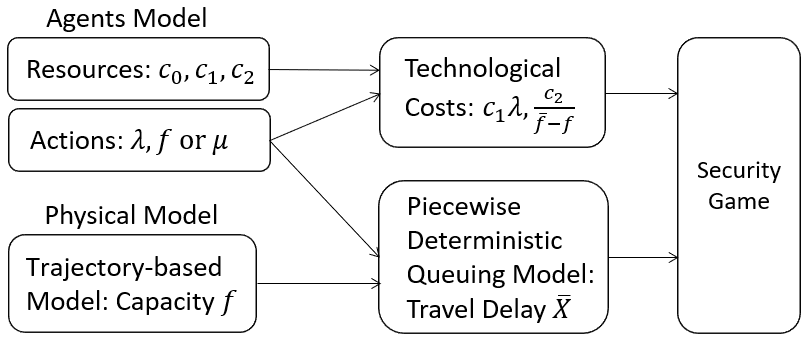}}
    \caption{The security game is based on multiple physical and agent models.}
    \label{fig:model}
\end{figure}

To address the above challenge, we develop a modeling approach that synthesizes three relevant models, viz., the trajectory-based model, the PDQ model, and security games, as shown in Fig.~\ref{fig:model}.
The trajectory-based model is typically used for vehicle coordination purpose, and the PDQ model is suitable for macroscopic performance analysis. 
We combine these models and strategic agent models to formulate a security game. The security game characterizes the interaction between the attacker and the SO.
The consequence of an attack is an increase in communication latency.
In our previous work~\cite{latency}, we developed a quantitative connection between intersection capacity and communication latency. Thus, the impact of an attack can be translated to a reduction in capacity. 
The intersection switches between the ``nominal mode" and the ``failure mode" in a Markovian manner, which turns out to be a PDQ model as studied in~\cite{PDQ}.

In the security game, the attacker is able to select the rate at which he/she jams the V2I communication. For practical purposes, we take the resource expenditure of attacks into consideration.
For the SO, we consider two types of actions: the SO is able to either select the saturation rate of the intersection, or the rate at which the compromised V2I communication is recovered. In the first case, the SO improves the saturation rate at the cost of increased investment in the V2I coordination hardware/infrastructure. In the second case, 
the SO improves the recovery rate at the cost of increased investment in the recovery capability, e.g., redundant sensing capabilities for attack detection and backup. 
The attacker (resp. the SO) is essentially interested in maximizing (resp. minimizing) the cost of queuing delay of the intersection, while reducing their respective technological investments.
For ease of presentation, we separately consider the two actions for the SO, which provides hints for the joint consideration of these two actions as future work.

Next, we study the structure of the equilibrium.
We analyze the security risk of the intersection by considering pure-strategy Nash Equilibria (NE) of the security game.
We first determine the regime where the Nash Equilibrium exists (Theorems 1 and 2). 
In practice, existence of a pure-strategy NE means that the SO can select a fixed level of investment to achieve security-by-design; otherwise, more sophisticated (e.g., feedback or randomized) strategies may be needed, which is a potential future work.
We find that if the technological cost for SO is relatively high, the NE exists.
When the NE exists, its structure is affected by the incoming traffic demand.
If the NE is non-trivial, the attacker will exhaustively use his/her budget, and the SO will either exhaustively use his/her budget or obtain the optimal action by applying the first-order optimality condition.
Moreover, if the attacker's budget is sufficiently high, then a pure-strategy NE does not exist. 
When the SO is able to select recovery rate, there exists another trivial NE: if the incoming traffic demand is small, the attacker does not attack and thus the SO does not recover. 
Our results provide hints for security risk analysis and help determine the optimal level of investment in vehicle
coordination hardware/infrastructure and recovery capabilities.

The main contributions of this paper are:
\begin{enumerate}
    \item 
    We provide a modeling approach to quantitatively analyze the security risk of signal-free intersections at the macroscopic level.
    \item 
    We formulate a security game in this scenario to characterize the interaction between the attacker and the SO.
    \item
    We study the equilibrium in the game to provide hints for preferable level of investment in vehicle coordination hardware/infrastructure and recovery capabilities.
    
\end{enumerate}

The rest of the paper is organized as follows. In Section II, we formulate the signal-free intersection scenario subject to jamming attacks, with the trajectory-based model and PDQ model introduced respectively. In Section III, the security game between the attacker and the SO is formulated, and the Nash Equilibrium in the security game is analyzed. In Section IV, we provide concluding remarks.

\section{Modeling Signal-free Intersections Subject to Jamming Attacks}
In this section, we consider of signal-free intersection based on V2I communications capacity, and discuss the possible jamming attacks that can cause latency. In Section A, we first introduce the trajectory-based model in the signal-free intersection. Then in Section B, we introduce the macroscopic PDQ model, and connect them together.
\subsection{Trajectory-based Model}
We consider a signal-free intersection with centralized vehicle coordination, where vehicles come along two orthogonal roads and need to pass through the intersection, as shown in Fig.~\ref{fig:intersec_fig}. A road-side unit (RSU) is established in the intersection to collect real-time kinematic information and coordinate the trajectories of traffic. 

In the signal-free intersection system, we consider a simplified trajectory coordination algorithm, which instructs each vehicle to track the pre-specified trajectory and to avoid collisions. It observes the actual location $s_i(t)$ and speed $v_i(t)$ for each vehicle $i$, and regulates the acceleration $u_i(t)$.
For each vehicle $i$ in the system, its kinematics evolves as follows:
\begin{align*}
s_{i}(t+1)&=s_{i}(t)+v_{i}(t) \delta,\\
v_{i}(t+1)&=v_{i}(t)+u_{i}(t) \delta+w_{i}(t),
\end{align*}
where $\delta$ is the time step size, and $w_{i}(t)$ is a disturbance term taking value in the range $[-\epsilon, \epsilon]$.
For all time $t$ and all vehicle $i$ on the same Origin-Destination (OD) pair, the safety constraint is imposed as follows:
$$
s_{i-1}(t)-s_{i}(t) \geq d+\beta v_{i}(t)
\quad\forall i.
$$
For vehicles on different ODs, they only need to satisfy the above when they are close to the crossing zone.
\begin{figure}[htbp]
    \centerline{\includegraphics[width=0.7\linewidth]{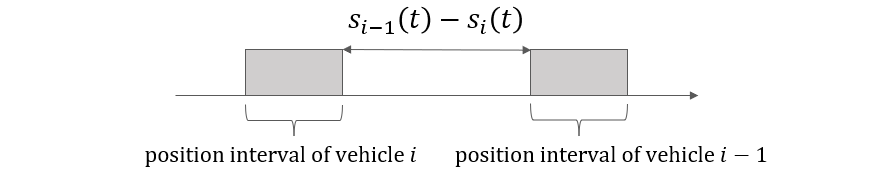}}
    \caption{The safety constraint of vehicles.}
    \label{fig:safety}
\end{figure}

Based on the above model, our previous work has deduced a simple disturbance-rejection controller in the context of signal-free intersections~\cite{latency}, which can be simplified as the following trajectory-tracking algorithm:
$$
u_{i}(t)=\frac{1}{\delta}\left(\frac{\bar{s}_{i}(t+2)-\left(s_{i}(t)+v_{i}(t) \delta\right)}{\delta}-v_{i}(t)\right),
$$
where $\bar{s}_{i}(t)$ is the reference location.
Then if we use $f$ to denote the capacity of the intersection with homogeneous vehicles, it can be given as follows:
$$
f=\frac{\bar{v}}{2 \epsilon \delta^{2}+d+\beta \bar{v}} \ [ veh /  sec ],
$$
where $\bar{v}$ is the nominal speed, $\epsilon$ is the maximal disturbance, $\delta$ is the coordination time step size, $d$ is the minimum static spacing that is considered to be safe between two vehicles, and $\beta$ is the reaction time.
Note that the SO can improve the saturation rate by increasing speed $\bar v$, reducing spacing $d$, and reducing coordination time step $\delta$; however, such improvement will require the SO to invest more on the hardware and infrastructure that support V2I connectivity.

\begin{figure}[htbp]
    \centerline{\includegraphics[width=0.5\linewidth]{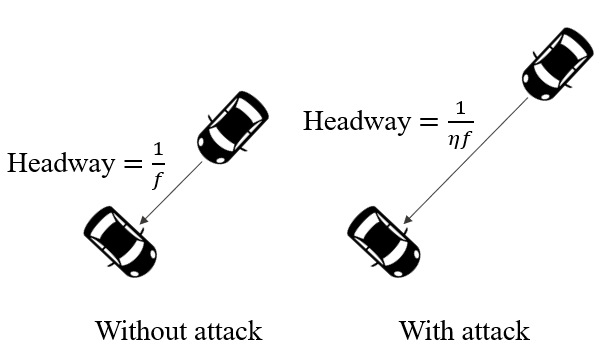}}
    \caption{The safe distance between vehicles increases when being attacked, resulting in a decrease in system capacity.}
    \label{fig:attack?}
\end{figure}

When malicious attacker jams the V2I communication channel, the SO has to adjust the coordination algorithm to accommodate to the degraded communication. For ease of presentation, we assume a simple reaction for the SO, viz. increasing the step size from $\delta$ to $\delta'$, without changing other parameters; see Fig.~\ref{fig:attack?}.
Hence, under jamming attacks, the saturation rate of the intersection changes to $f'$ as follows:
$$
f^{\prime}=\frac{\bar{v}}{2 \epsilon(\delta')^2+d+\beta \bar{v}} \ [ veh /  sec ].
$$
We also denote the ratio of the capacity with attack to the capacity without attack as a discounting ratio $\eta = \frac{f'}{f}$, which is given by
\begin{align*}
    \eta=\frac{2\epsilon\delta^2+d+\beta\bar v}{2\epsilon(\delta')^2+d+\beta\bar v}
    =1-\frac{2\epsilon((\delta')^2-\delta^2)}{2\epsilon(\delta')^2+d+\beta\bar v}.
\end{align*}
We assume $\eta$ to be a constant in this paper, but one can indeed generalize the formulation to consider $\eta$ as a decision variable.

\subsection{Piecewise Deterministic Model (PDQ)}
\label{sub_pdq}

\begin{figure}[htbp]
    \centerline{\includegraphics[width=0.4\linewidth]{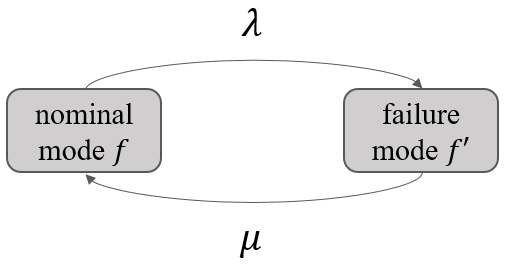}}
    \caption{The two-state Markov chain.}
    \label{fig:markov}
\end{figure}

The two different modes of capacity in the system can be represented as a two-state Markov chain, as shown in Fig.~\ref{fig:markov}. We denote that one state is the ``nominal mode" of the intersection, with capacity $f$, and the other state is the ``failure mode" under attack, with capacity $f'$. We define a discrete state $Y \in\{0,1\}$ to represent the current state, with 0 as the ``nominal mode" and 1 as the ``failure mode". In the two-state Markov chain, $\lambda$ represents the transition rate from ``nominal mode" to ``failure mode", and $\mu$ represents the transition rate from ``failure mode" to ``nominal mode". Then the steady-state possibilities of the game is given as $\pi=\left[\frac{\mu}{\lambda+\mu}, \frac{\lambda}{\lambda+\mu}\right]$.

\begin{figure}[htbp]
    \centerline{\includegraphics[width=0.5\linewidth]{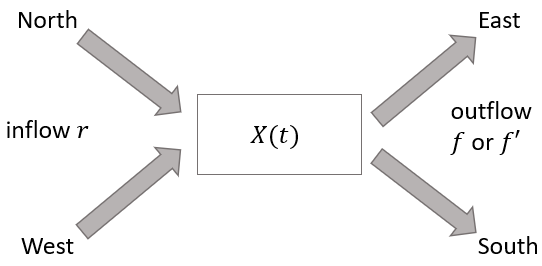}}
    \caption{Illustration of inflow and outflow.}
    \label{fig:f&r}
\end{figure}

In a PDQ model, consider a simple scenario shown in Fig.~\ref{fig:f&r}~\cite{PDQ}. The inflow $r$ represents the total inflow of the system, which is the sum of traffic flows coming north and west in $[ veh /  sec ]$. Similarly, the outflow $f$ or $f'$ represents the total outflow, which is sum of traffic flows heading to east and south in $[ veh /  sec ]$.

Then denote the traffic queue size in the system as continuous state $x \in \mathbb{R}_{\geq 0}$, the system dynamics can be described as follows:
$$
\frac{d}{d t} X(t)= \begin{cases}r-f<0 & Y(t)=0, 
\\ r-f^{\prime}>0 & Y(t)=1.\end{cases}
$$
\begin{figure}[htbp]
    \centerline{\includegraphics[width=0.7\linewidth]{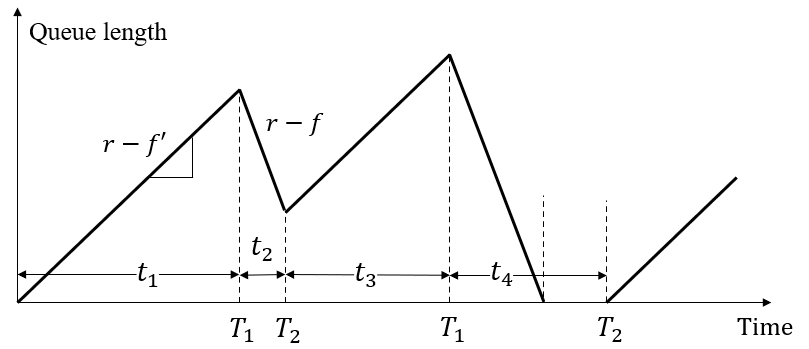}}
    \caption{Illustration of the piecewise-deterministic queuing model.}
    \label{fig:PDQ}
\end{figure}

Consider the long-time average queue length in the system, which can reflect the capacity and performance of signal-free intersections. Intuitively, when the system is under attack in the ``failure mode", and inflow $r$ is larger than current outflow $f'$, the queue length keeps increasing. On the contrary, when the system is in the ``nominal mode", and the inflow $r$ is smaller than current outflow $f$, the queue lengths keeps decreasing till zero. Illustration of a simple case is shown in Fig.~\ref{fig:PDQ}. We can see that in time period $t_1$ and $t_3$, the queue length increases at the rate of $r-f'$. At time $T_1$, it reaches the turning point, and then decreases at the rate of $r-f$ in the time period $t_2$ and $t_4$, respectively.

With $\lambda$ and $\mu$ being the given transition rates, the long-time average queue length can be developed to quantify the queuing delay cost at the intersection as follows:
$$ \bar{X}(\lambda, \mu, f) =\left\{\begin{array}{cc}
     0 & 0 \leq r \leq f^{\prime}, \\
    \frac{\lambda}{\mu} \frac{\left(r-f^{\prime}\right)\left(f^{\prime}-f\right)}{\lambda\left(r-f^{\prime}\right)+\mu(r-f)} & f^{\prime}<r<\frac{\mu}{\lambda+\mu} f+\frac{\lambda}{\lambda+\mu} f^{\prime}, \\
\infty & o . w ..
\end{array}\right.
$$
For the derivation of the above formula, see \cite{kulkarni97}.

\section{Security Game for Jamming Attacks}
In this section, we formulate a security game between the attacker and the SO in the system, analyze the equilibrium of the game, and derive practical insights. In Section A, we define the agent models and their utilities. In Section B, we study the case when the SO is able to select the saturation rate of the intersection. Then in Section C, we study the case when the SO is able to select the recovery rate.

\subsection{Security Game}

We consider that the attacker is able to select the Poisson rate $\lambda$ of jamming attack to transit the ``nominal mode" into the ``failure mode" in a given range [0, $\bar{\lambda}$]. That is, the time interval between the last attack getting resolved and the next attack occurring is a random variable following the exponential distribution with parameter $\lambda.$ Hence, one can interpret $\lambda$ as the ``attacking rate'' of the attacker, and the upper bound $\bar{\lambda}$ as the technological limit of the attacking rate.
In addition, typical anti-jamming methods such as channel hopping and spread spectrum techniques are designed to increase the attacker's jamming power~\cite{hopping1},~\cite{hopping2},~\cite{anti-jam},~\cite{spread}. 
Hence, in order to generate more rational and practical decisions, we take the resource expenditure of jamming attacks (as well as SO decisions) into consideration. 
Specifically, we assume that every attack will cost the attacker $c_1$ amount of money. One can indeed also consider $\delta'$, i.e. the degraded time step size as a decision variable for the attacker; this can be interpreted as the ``intensity'' of each attack. We do not consider this scenario in this paper, but our results can be extended to cover it.

For the SO, we consider two actions separately. In the first scenario, the SO is able to select the saturation rate of the intersection $f$ within range [0, $\bar{f}$], where $\bar{f}$ is the technological limit of the saturation rate of the intersection. 
As discussed in Section~\ref{sub_pdq}, the SO can improve $f$ by investing in finer V2I coordination hardware/infrastructure that enables higher crossing speed and/or shorter inter-vehicle headways. One can view the selection of $f$ as a collective selection of the above-mentioned lower-level specifications.
The selection of $f$ is indeed with a technological cost; we consider this cost to be $c_2/(\bar f-f)$, which means that approaching the technological limit $\bar f$ will be extremely costly.
In addition, we assume that the recovery rate takes a nominal value $\mu_0$. Thus, the utility for the attacker $u_1(\lambda, f)$ is the summation of two parts: (i) the queuing cost $c_{0}\bar{X}(\lambda,\mu_0, f)$, which is the product of the long-time average queuing delay at the intersection and the value of time, (ii) $-c_1\lambda$, which is the expectation of the attacking cost for the attacker. The utility for the SO $u_2(\lambda, f)$ is also the summation of two parts: (i) the inverse of the queuing cost $-c_{0}\bar{X}(\lambda, \mu_0, f)$, and (ii) $-\frac{c_2}{\bar{f}-f}$, which captures the technological cost. Hence, the utility functions of the attacker and the SO are given as follows:
\begin{align*}
\mbox{Attacker: }u_{1}(\lambda, f)&=c_{0}\bar{X}(\lambda, \mu_0, f)-c_{1} \lambda,\\
\mbox{SO: }u_{2}(\lambda, f)&=-c_{0}\bar{X}(\lambda,\mu_0, f)-\frac{c_2}{\bar{f}-f}.
\end{align*}

In the second scenario, the SO is able to select the Poisson rate $\mu$ of recovery to transit the ``failure mode" into the ``normal mode" in a given range [0, $\bar{\mu}$]. That is, the duration of each attack is a random variable following the exponential distribution with parameter $\mu$. The upper bound $\bar\mu$ can be interpreted as the technological limit of any recovery capability. The magnitude of $\mu$ can be interpreted as the SO's investment in the capability of recovering the intersection from an attack; examples of such capabilities include redundant sensing capabilities for attack detection and backup, more resilient V2I connectivity that can quickly and/or automatically fix jamming, and human labor for inspection and restoration.
When the SO selects $\mu$, we assume that the saturation rate takes a nominal value $f_0$. The utility functions for the attacker and the SO are thus given by
\begin{align*}
\mbox{Attacker: }u_{1}'(\lambda, \mu)&=c_{0}\bar{X}(\lambda, \mu, f_0)-c_{1} \lambda,\\
\mbox{SO: }u_{2}'(\lambda, \mu)&=-c_{0}\bar{X}(\lambda, \mu, f_0)-c_{2}' \mu.
\end{align*}

With both the attacker and the SO aiming to maximize their own utility function, a pure-strategy Nash Equilibrium (NE) may exist in the security game. Recall that a pair ($\lambda^*$, $f^*$) $\in [0, \bar{\lambda}] \times[0, \bar{f}]$ (resp. ($\lambda^*$, $\mu^*$) $\in [0, \bar{\lambda}] \times[0, \bar{\mu}]$) is said to be an NE if in ($\lambda^*$, $f^*$) (resp. ($\lambda^*$, $\mu^*$)), neither the attacker nor the SO has an incentive to unilaterally deviate from their chosen strategy after considering an opponent's choice, where $\lambda^*$ and $f^*$ (resp. $\mu^*$) is the choice of the attacker and the SO at equilibrium, respectively.

One can indeed allow the SO to select $\mu$ and $f$ simultaneously. However, such a scenario will lead to more involved math. We do not consider this complexity in this paper, but the results in the following provide hints for that more complex scenario.

\subsection{Equilibrium When the SO Selects $f$}
In this section, we consider the case when the SO is able to select the saturation flow rate of the intersection $f$. We provide an analytical characterization of the Nash Equilibrium of the security game as described above. Next, we focus on analyzing the equilibrium regimes for various costs $c_1$ and $c_2$. 

\textit{Theorem 1:} 
With $\mu$ fixed, a Nash Equilibrium $\left(\lambda^{*}, f^{*}\right)$ exists if and only if both conditions (1a) and (1b) are satisfied:
\begin{subequations}
    \begin{align}
        &\bar{\lambda}<\frac{\mu(r-\bar{f})}{-\left(r-\eta\bar{f}\right)},    \tag{1a}\\
        &c_{1} \leq  \frac{c_0(\eta\widehat{f}-\widehat{f})(r-\eta\widehat{f})}{\mu\left[\bar{\lambda}\left(r-\eta\widehat{f}\right)+\mu(r-\widehat{f})\right]}, \tag{1b}   
    \end{align}
\end{subequations}

where $\widehat{f}$ denotes the solution of the equation $\frac{\partial u_{2}(\bar{\lambda}, f)}{\partial f}=0$ that satisfies $f \in(\frac{\bar{\lambda}+\mu}{\eta\bar{\lambda}+\mu}r, \bar{f}]$. 
Furthermore, if an NE exists, then it is unique and given by:
$$
\left(\lambda^{*}, f^{*}\right)= (\bar{\lambda}, \widehat{f}).
$$

\begin{figure}[htbp]
        \centerline{\includegraphics[width=0.5\linewidth]{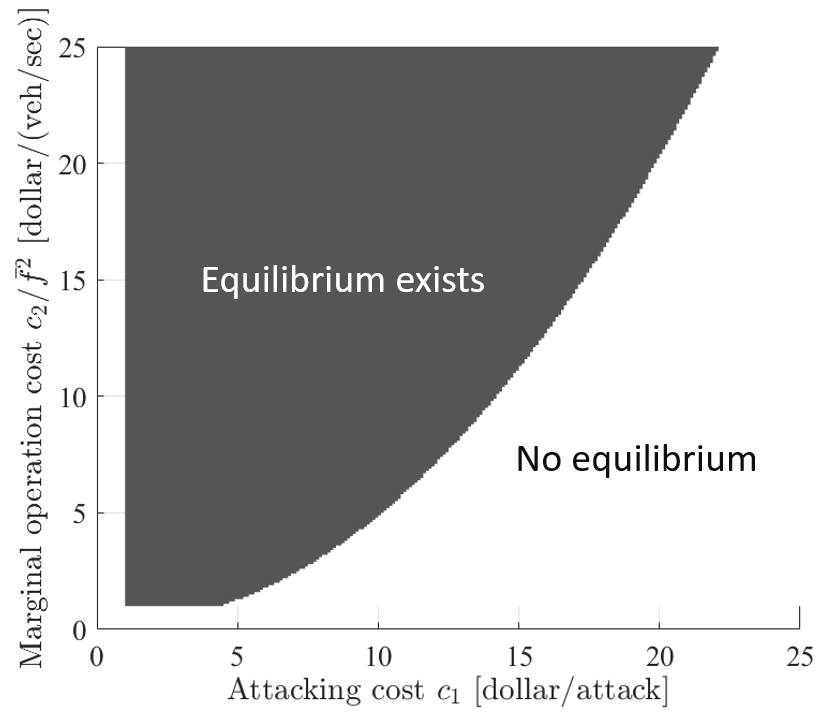}}
    \caption{Regimes for the existence of pure-strategy NE when the SO selects $f$.}
    \label{fig:regime c1}
\end{figure}

Using Theorem 1, we can analyze the regimes for the technological costs $c_1$, $c_2$ under which an NE exists. Consider a typical example with parameters given in Table~\ref{tab1}. Here we take the nominal recovery rate $\mu_0$ into consideration. The regime diagram is shown in Fig.~\ref{fig:regime c1}.

\begin{table}[h!]
\centering
\begin{tabular}{cccc} 
 Parameter Name & Notation & Value & Unit\\
 \hline
 value of time & $c_0$ & 1 & dollar/sec \\
 maximal attacking rate & $\bar{\lambda}$ & 0.1 & per sec\\
 discounting ratio & $\eta$ & 0.4 & N/A\\
 maximal saturation rate & $\bar{f}$ & 1 & veh/sec\\
 inflow traffic rate & $r$ & 0.3 & veh/sec\\
 nominal recovery rate & $\mu_0$ & 0.2 & per sec\\ 
 maximal recovery rate & $\bar{\mu}$ &0.4 & per sec\\
 nominal saturation rate & $f_0$ & 0.5 & veh/sec\\
 \hline
\end{tabular}
\caption{Parameters in the example.}
\label{tab1}
\end{table}
According to observations of Fig.~\ref{fig:regime c1}, we can generate the following implications of the existence of NE: 

\begin{enumerate}
    \item 
    When the operation cost $c_2$ is fixed, ($c_2$ is the marginal operation cost when $f=0$) and the attacking cost $c_1$ is low, the equilibrium exists. Intuitively, the low cost of attack results in the high utility value of the attacker, hence the attacker has a relatively large incentive to attack, the equilibrium exists.
    \item 
    When $c_2$ is fixed and $c_1$ is high, the equilibrium does not exist. Intuitively, the high cost of attack results in the low utility value of the attacker, hence the attacker does not have the incentive to attack. Then the queuing cost would be zero, in order to minimize technological costs, the SO would not let any vehicle pass the intersection. However, this motivates the attacker to do some attacks, which causes the SO to choose a higher saturation flow rate. 
    Hence, in this case, a static equilibrium can not be reached.
    \item 
    When $c_1$ is fixed and $c_2$ is low, the equilibrium does not exist. Intuitively, the low operating cost enables the SO to always afford the cost and keep the system resilient to attacks. Therefore, each time the attacker changes the attacks, the SO will change its saturation rate selection, the equilibrium exists.
    \item 
    When $c_1$ is fixed and $c_2$ is high, the equilibrium exists. Intuitively, the high operation cost results in the low utility value of the SO, then the SO is not always able to afford the cost, and would choose to tolerate some of the attacks. Therefore, a static equilibrium can be reached.
\end{enumerate}

\textit{Proof of Theorem 1.}
We first prove the \textbf{sufficiency} of Theorem 1. 
Assume that both conditions (1a) and (1b) are satisfied, then we respectively analyze the utility functions for the attacker and the SO to derive the unique NE $(\bar{\lambda}, \widehat{f})$. 
When condition (1a) is satisfied, the long-time average queue length $\bar{X}(\lambda, \mu, f)$ is not infinity. If $\bar{X}(\lambda, \mu, f) = 0$, then the maximum utility for both the attacker and the SO can be achieved when $\lambda = f = 0$. However, when $f=0$, the intersection no longer has practical significance. Hence, when $\bar{X}(\lambda, \mu, f) = 0$, the equilibrium does not exist. 
Then we consider the case when $\bar{X}(\lambda, \mu, f)$ is neither infinity nor zero.
Take the first and second partial derivatives of $u_1(\lambda, f)$ with respect to $\lambda$, we can find that $\frac{\partial \bar{X}(\lambda, \mu, f)}{\partial \lambda} > 0$, $\frac{\partial u_{1}(\lambda, f)}{\partial \lambda} = c_0\frac{\partial \bar{X}(\lambda, \mu, f)}{\partial \lambda} - c_1$,
$\frac{\partial^2 u_{1}(\lambda, f)}{\partial \lambda^2} = c_0\frac{\partial^2 \bar{X}(\lambda, \mu, f)}{\partial \lambda^2} > 0$, hence $u_1(\lambda, f)$ is convex. 
To maximize the utility, we only need to compare the values on the boundaries, which are $u_1(0, f)$ and $u_1(\bar{\lambda}, f)$. When condition (1b) is satisfied, $u_1(\bar{\lambda}, f) > u_1(0, f)$ is always true for all $f \in [0, \bar{f}]$. Hence, the maximum utility for the attacker is achieved when $\lambda = \bar{\lambda}$, therefore $\lambda^* = \bar{\lambda}$. 
Then we analyze the utility function for the SO $u_2(\lambda, f)$. When $\bar{X}(\lambda, \mu, f)$ is neither infinity nor zero, $f$ is within the range $(\frac{\bar{\lambda}+\mu}{\eta\bar{\lambda}+\mu}r, \bar{f}]$.
Similarly, take the first and second partial derivatives of $u_2(\lambda, f)$ with respect to $f$, we can find that 
$\frac{\partial \bar{X}(\lambda, \mu, f)}{\partial f} < 0$,
$\frac{\partial u_2(\lambda, f)}{\partial f} =  -c_0\frac{\partial \bar{X}(\lambda, \mu, f)}{\partial f}- \frac{c_2}{(\bar{f}-f)^2}$, $\frac{\partial^2 u_{2}(\lambda, f)}{\partial f^2} < 0$, therefore the utility function $u_2(\lambda, f)$ is concave. 
To maximize the utility and find the equilibrium, we want to solve $\frac{\partial u_{2}(\bar{\lambda}, f)}{\partial f}=0$. Denote the left side of the equation after removing the denominators as $g(f)$, we can find that $g(\frac{\bar{\lambda}+\mu}{\eta\bar{\lambda}+\mu}r) > 0$ and $g(\bar{f}) < 0$. Hence, the unique solution of the above function always exists within the range $(\frac{\bar{\lambda}+\mu}{\eta\bar{\lambda}+\mu}r, \bar{f}]$, which indicates the existence of the equilibrium. With the solution denoted as $\widehat{f}$, $f^* = \widehat{f}$ would be the optimal choice for the SO. 
Hence, when conditions (1a) and (1b) are satisfied, there exists an unique NE $(\bar{\lambda}, \widehat{f})$. The sufficiency of Theorem 1 has been proved.

Next, we prove the \textbf{necessity} of Theorem 1 by contradiction.
If the condition (1a) is not satisfied, assume that the NE $\left(\lambda^{*}, f^{*}\right)$ exists, then $r > \frac{\mu}{\bar{\lambda}+\mu} \bar{f}+\frac{\bar{\lambda}}{\bar{\lambda}+\mu} \eta \bar{f} \geq \frac{\mu}{\bar{\lambda}+\mu} f^{*}+\frac{\bar{\lambda}}{\bar{\lambda}+\mu} \eta f^{*}$. Hence, $\bar{X}\left(\lambda^{*},\mu, f^{*}\right) = \infty$, the choices of the attacker and the SO can not significantly influence the utility, the NE does not hold. 
If the (1b) is not satisfied, assume that the NE $\left(\lambda^{*}, f^{*}\right)$ exists. Since the optimal choice for the attacker is either $0$ or $\bar{\lambda}$, $u_1\left(\lambda^{*}, f^{*}\right) = c_{0}\bar{X}(\bar{\lambda}, \mu, f^{*}) - c_1\bar{\lambda} < 0 = u_1\left(0, f^{*}\right)$, therefore $\lambda^{*} = 0$. However, as stated above, if $\lambda^{*} = 0$, then $f^{*} = 0$, the intersection no longer has practical significance, which leads to a contradiction.
Hence, the NE exists only if both conditions (1a) and (1b) are satisfied. The necessity of Theorem 1 has been proved.
In conclusion, the unique NE $\left(\lambda^{*}, f^{*}\right)= (\bar{\lambda}, \widehat{f})$ exists if and only if both conditions (1a) and (1b) are satisfied.{\hfill $\square$\par}

\subsection{Equilibrium When the SO Selects $\mu$}
In this section, we consider the case when the SO is able to select the recovery rate $\mu$, analyze the structure of the equilibrium and provide hints.

\textit{Theorem 2:} 
With $f$ fixed, a Nash Equilibrium  $\left(\lambda^{*}, \mu^{*}\right)$ exists if and only if one of the following conditions hold:
\begin{enumerate}
    \item Underloaded system:
    \begin{subequations}
    \begin{align}
        & 0 \leq r \leq f^{\prime}; \tag{2.1}
    \end{align}
    \end{subequations}
    
    \item Redundant recovery resource: 
    \begin{subequations}
    \begin{align}
        &\bar{\lambda}<\frac{\bar{\mu}(r-f)}{-\left(r-f^{\prime}\right)}, \tag{2.2a} \\
        &\bar{\mu}\geq \widehat{\mu}, \tag{2.2b}\\
        &c_{1} \leq  \frac{c_0\left(r-f^{\prime}\right)\left(f^{\prime}-f\right)}{\widehat{\mu}	\left[\bar{\lambda}\left(r-f^{\prime}\right)+\widehat{\mu}(r-f)\right]}; \tag{2.2c}
    \end{align}
    \end{subequations}
    \item Exhaustive recovery resource: 
    \begin{subequations}
    \begin{align}
        &\bar{\lambda}<\frac{\bar{\mu}(r-f)}{-\left(r-f^{\prime}\right)}, \tag{2.3a} \\
        &\bar{\mu}<\widehat{\mu}, \tag{2.3b}\\
        &c_{1} \leq  \frac{c_0\left(r-f^{\prime}\right)\left(f^{\prime}-f\right)}{\bar{\mu}	\left[\bar{\lambda}\left(r-f^{\prime}\right)+\bar{\mu}(r-f)\right]}, \tag{2.3c}
    \end{align}
    \end{subequations}
\end{enumerate}
where $\widehat{\mu}$ denotes the solution of the equation $\frac{\partial u_{2}(\bar{\lambda}, \mu)}{\partial \mu}=0$ that satisfies $\mu \in[0, \bar{\mu}]$. 
Furthermore, if an NE exists, then it is unique and given by:
$$
\left(\lambda^{*}, \mu^{*}\right)= 
\begin{cases}(0,0) & \text{condition}(2.1), \\ 
(\bar{\lambda}, \widehat{\mu}) & \text{condition}(2.2) , \\
(\bar{\lambda}, \bar{\mu}) & \text{condition}(2.3).
\end{cases}
$$

\begin{figure}[htbp]
        \centerline{\includegraphics[width=0.5\linewidth]{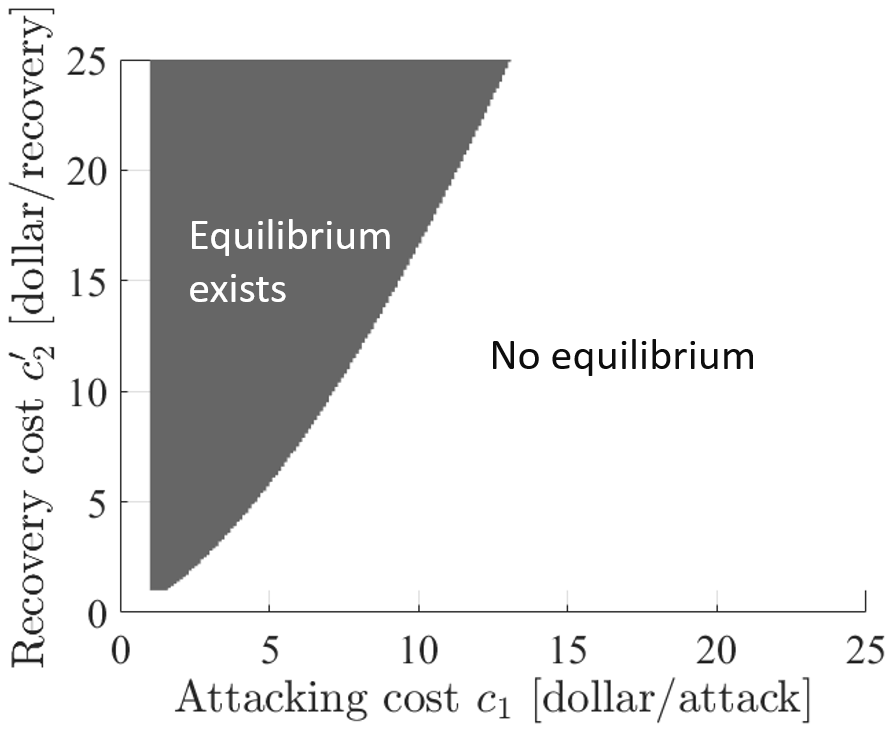}}
    \caption{Regimes for the existence of pure-strategy NE when the SO selects  $\mu$.}
    \label{fig:regime c2}
\end{figure}

Similarly, using Theorem 2, we can analyze the regimes for the technological costs $c_1$, $c_2'$ under which an NE exists. Consider the same example with parameters given in Table~\ref{tab1}. Here we take the nominal saturation rate $f_0$ into consideration. The regime diagram is shown in Fig.~\ref{fig:regime c2}, which implies the following.

\begin{enumerate}
    \item 
    When the recovery cost $c_2'$ is fixed and the attacking cost $c_1$ is low, the attacker has a relatively large incentive to attack, the equilibrium exists.
    \item 
    When $c_2'$ is fixed and $c_1$ is high, the attacker does not have the incentive to attack, so the SO does not need to conduct recovery operations, which motivates the attacker to do some attacks, however. Hence, a static equilibrium can not be reached.
    \item 
    When $c_1$ is fixed and $c_2'$ is low, the low recovery cost enables the SO to always afford the recovery cost and change its recovery rate selection whenever the attacker changes the attacks, so a static equilibrium can not be reached.
    \item 
    When $c_1$ is fixed and  $c_2'$ is high, the SO is not always able to afford the cost, and would choose to tolerate some of the attacks, a static equilibrium can be reached.
\end{enumerate}

\textit{Proof of Theorem 2.}
We first prove the \textbf{trivial case} when the NE is (0, 0). Then we jointly consider the other two cases. 
Firstly, in the underloaded system, if the condition (2.1) is satisfied, the long-time average queue length $\bar{X}(\lambda, \mu, f) = 0$. To maximize the utility functions, both the attacker and the SO have no incentive to work. Hence, in this case, the unique NE ($\lambda^*$, $\mu^*$) = (0, 0) exists.

Secondly, we jointly consider the rest two cases and prove the \textbf{sufficiency}. Assume that both conditions (2.2) and (2.3) are satisfied, then we respectively analyze the utility functions to derive the unique NE. 
When condition (2.2a) (or condition (2.3a)) is satisfied, $\bar{X}(\lambda, \mu, f)$ is neither zero nor infinity.
Similarly, take the first and second partial derivatives of $u_1'(\lambda, \mu)$ with respect to $\lambda$, we can find that $\frac{\partial \bar{X}(\lambda, \mu, f)}{\partial \lambda} > 0$, $\frac{\partial u_{1}'(\lambda, \mu))}{\partial \lambda} = c_0\frac{\partial \bar{X}(\lambda, \mu, f)}{\partial \lambda} - c_1$,
$\frac{\partial^2 u_{1}'(\lambda, \mu)}{\partial \lambda^2} > 0$,  $u_1'(\lambda, \mu))$ is convex. The maximum utility is either at $\lambda = \bar{\lambda}$ or $\lambda = 0$. 
If $\lambda = 0$, $\bar{X}(\lambda, \mu, f) = 0$, in order to maximize the utility $u_2'(\lambda, \mu)$, the SO should choose $\mu = 0$. However, such circumstance does not lie in the range of the second case of $\bar{X}(\lambda, \mu, f)$, which leads to contradiction.
Hence, if the NE exists, $\lambda^{*} = \bar{\lambda}$.
Take the first and second partial derivatives of $u_2'(\lambda, \mu)$ with respect to $\mu$, we can find that 
$\frac{\partial \bar{X}(\lambda, \mu, f)}{\partial \mu} < 0$,
$\frac{\partial u_2'(\lambda, \mu)}{\partial \mu} =  -c_0\frac{\partial \bar{X}(\lambda, \mu, f)}{\partial \mu}- c_2'$, $\frac{\partial^2 u_{2}'(\lambda, \mu)}{\partial \mu^2} < 0$, therefore the utility function $u_2'(\lambda, \mu)$ is concave. 
To maximize the utility and find the equilibrium, we want to solve $\frac{\partial u_{2}'(\bar{\lambda}, \mu)}{\partial \mu}=0$. Similarly, since the above equation has and only has one solution within the range $(\frac{\bar{\lambda}(r-f')}{-(r-f)}, +\infty]$, we denote the unique solution as $\widehat{\mu}$. 
According to the characteristics of the function $u_2'(\lambda, \mu)$, when condition (2.2b) is satisfied, the maximum can be achieved at $\mu = \widehat{\mu}$, otherwise condition (2.3b) is satisfied, and the maximum can be achieved at $\mu = \bar{\mu}$. 
Consider the utility for the attacker, when both conditions (2.2b)  and (2.2c) are satisfied, or when conditions (2.3b)  and (2.3c) are satisfied, $u_1'(\bar{\lambda}, \mu) > u_1'(0, \mu)$ is always true for any $\mu in [0, \bar{\mu}]$.
Hence, when conditions (2.2) are satisfied, there exists an unique NE $(\bar{\lambda}, \widehat{\mu})$, and when conditions (2.3) are satisfied, there exists an unique NE $(\bar{\lambda}, \bar{\mu})$. The sufficiency has been proved.

Next, we prove the \textbf{necessity} of the two cases by contradiction.
If the condition (2.2a) (or condition (2.3a)) is not satisfied, assume that the NE $\left(\lambda^{*}, f^{*}\right)$ exists, similarly, $\bar{X}\left(\lambda^{*},\mu, f^{*}\right) = \infty$, the NE does not hold. 
Conditions (2.2b) and (2.3b) are complementary to each other, so either of it has to be satisfied.
If condition (2.2b) holds and (2.2c) is not satisfied, assume that the NE exists, then $u_1'\left(\lambda^{*}, \mu^{*}\right) = c_{0}\bar{X}(\bar{\lambda}, \mu, f^{*}) - c_1\bar{\lambda} < 0 = u_1'\left(0, \mu^{*}\right)$, which leads to a contradiction. Condition (2.3c) can be proved in the same way. Hence, the NE exists only if conditions (2.2) or (2.3) are satisfied. The necessity  has been proved.
To summarize, the Theorem 2 has been proved.{\hfill $\square$\par}

\section{Concluding Remarks}

In this paper, we provide a modeling approach to model the impact of jamming attacks on the intersection-level performances for connected and autonomous vehicles in signal-free intersections. Based on the model, we formulate a security game to analyze the interaction between the attacker and the system operator. By studying the equilibrium structure, we provide hints for security risk analysis and resource allocation for vehicle
coordination hardware/infrastructure and recovery capabilities.

Future researches can be conducted to generalize our results under multiple contexts. One possible direction is to consider an incomplete information security game. Our work considers a complete information game between the attacker and the SO, however, in reality it is often difficult to obtain complete real-time information. Furthermore, the sequential security game between the attacker and the SO, and the scenario when the SO is able to simultaneously select both parameters $\mu$ and $f$ are also worth studying.

\bibliographystyle{IEEEtran}
\bibliography{main}










\end{document}